\documentclass[11pt]{article}
\title {Torus Knot complements:\\ 
A natural series for the natural logarithm}
\author
{Oliver T. Dasbach
\thanks{e-mail: kasten@math.lsu.edu, 
http://www.math.lsu.edu/$\sim$\!\! kasten}
\thanks{partially supported by NSF DMS-0306774  and DMS-FRG-0456275 }
\\Louisiana State University\\
Department of Mathematics\\Baton Rouge, LA 70803
}

\date{}

\usepackage{amsmath,amssymb,amsthm, graphics} 

\newtheorem{theorem}{Theorem}[section]

\newtheorem{lemma}[theorem]{Lemma}

\newtheorem{remark}[theorem]{Remark}

\newcommand{\C} {\mathbb{C}} 
\newcommand{\Z} {\mathbb{Z}}

\newcommand{\isom}{\cong}

\newcommand{\tr} {\mbox{\rm tr}}

\newcommand{\Vol} {\mbox{\rm Vol}}
\newcommand{\Det} {\mbox{\rm det}}

\begin{document} 
\maketitle
\begin{abstract} 
L\"uck expressed the Gromov norm of a knot complement in terms of an infinite series that can be
computed from a presentation of the fundamental group of the knot complement.

In this note we show that L\"uck's formula, applied to torus knots, yields surprising power series expansions 
for the logarithm function. This generalizes an infinite series of Lehmer for the natural logarithm of $4$.
\end {abstract}

\section{Background}

First, we need to fix some notations. 
Let $K$ be a knot and $$G=\pi_1(S^3-K)=\langle x_1, \dots, x_g | r_1, \dots, r_{g-1} \rangle$$
be a presentation of the fundamental group of the knot complement.
For a square matrix $M$ with entries in $\C G$ the trace $\tr_{\C G}(M)$
denotes the coefficient of the unit element in the sum of the diagonal elements. 
The matrix $A^*=(\bar{a}_{j,i})$ is the conjugate transpose of $A=(a_{i,j})$ with 
conjugation
$$\overline {\sum_{g \in G} c_g g} = \sum_{g \in G} \overline {c_g} g^{-1}.$$   

Let $$F=\left ( 
\begin{array}[pos]{ccc}
	\frac{\partial r_1} {\partial x_1} & \dots & \frac{\partial r_1} {\partial x_g}\\
	\vdots &\ddots & \vdots\\
	\frac{\partial r_{g-1}} {\partial x_1} & \dots & \frac{\partial r_{g-1}} {\partial x_g}\\
\end{array}
\right )$$
be the Fox Jacobian (e.g. \cite{BZ}) of the presentation.
We obtain a $(g-1) \times (g-1)$-matrix $A$ by deleting one of the columns of $F$.

\begin{theorem}[L\"uck \cite{LueckBook:L2invariants}] \label{LuecksTheorem} Suppose the deleted column in the Fox Jacobian does not correspond to
a generator $x_i$ of $G$ that represents a trivial element in $G$.  Then for a hyperbolic knot $K$ and for 
$k$ sufficiently large
it holds:

\begin{equation} \label{LuecksFormula}
\frac 1 {3 \pi} \Vol(S^3-K)= 2 (g-1) \ln (k) - \sum_{n=1}^{\infty} \frac 1 n \tr_{\C G} \left ( ( 1-k^{-2} A A^*)^n \right ) 
\end{equation}
\end{theorem}

The value of $k$ can be chosen to be the product of $(g-1)^2$ and the maximum of the $1$-norm of the entries in $A$ (see \cite{Lueck:Proceedings}). 

In case of a non-hyperbolic knot, the right hand side of (\ref{LuecksFormula}) is proportional to the Gromov norm of the knot complement.

To give a hands-on example for L\"ucks formula, we elementary compute the case where $A$ is a matrix with complex entries: 

\begin{lemma} \label{complex entries}
Suppose $A=(g-1)\times(g-1)$ is non-singular with entries in $\C$.

If $k \geq ||A||_2$, i.e. $k^2$ is greater or equal to the largest eigenvalue of $A A^*$, then
\begin{eqnarray*}
 2 (g-1) \ln (k) 
- \sum_{n=1}^{\infty} \frac 1 n \tr \left 
( ( I-k^{-2} A A^* )^n \right ) &=& 2 \ln | \det A |
\end{eqnarray*}

\end{lemma}

\begin{proof}
Let $\lambda_1, \dots, \lambda_{g-1}$ be the eigenvalues of $A A^*$.
Since $A A^*$ is diagonizable, we have:
\begin{eqnarray*}
&&(g-1) \ln k^2 - \sum_{n=1}^{\infty} \frac 1 n \tr (I - k^{-2} A A^*)^n\\
&=& (g-1) \ln k^2 - \sum_{n=1}^{\infty} \frac 1 n \sum_{j=1}^{g-1} 
\left (1- \frac {\lambda_j}{k^2} \right )^n\\
&=& (g-1) \ln k^2 - \sum_{j=1} ^{g-1} \sum_{n=1}^{\infty} 
\frac 1 n \left (1- \frac {\lambda_j}{k^2} \right )^n\\
&=& (g-1) \ln k^2 + \sum_{j=1}^{g-1} \ln \left (\frac {\lambda_j}{k^2} \right )\\
&=& \ln \det (A A^*)\\
&=& 2 \ln |\det A|
\end{eqnarray*}

\end{proof}

Lemma \ref{complex entries} has a generalization to matrices over a polynomial ring in one variable:

\begin{theorem} [\cite{LueckBook:L2invariants}, see also \cite{Deninger:FugledeKadison}] \label{MahlerMeasure}
Let $A$ have entries in $\C[x_1, \dots x_s]$. Then the right hand side of (\ref{LuecksFormula})
equals:

$$2 m(\Det(A)),$$
where $m(p)= \int_0^1 \dots \int_0^1 \ln |p(e^{2 \pi i t_1}, \dots, e^{2 \pi i t_s})| dt_1 \dots dt_s$ is the (logarithmic) Mahler measure (see e.g. \cite{EverestWard:MahlerMeasure})
of the polynomial $p(x_1, \dots, x_s)$.
\end{theorem}

\section{A power series for the natural logarithm}

Let $K=T(p,q)$ be the $(p,q)$-torus knot. It is well-known (e.g. \cite{BZ}) that a presentation for the knot group of $T(p,q)$ is given by:

$$G=\pi_1(S^3-T(p,q))=\langle a, b| a^p=b^q \rangle.$$

Hence, the Fox Jacobian for $G$ is
$$F=(1+a+a^2+ \dots + a^{p-1}, -1 - b - \dots - b^{q-1})$$  

and we can chose $A$ to be the $1 \times 1$ matrix:
$$A=(1+a+\dots+a^{p-1}) \qquad \mbox{ and thus } \qquad A^*=(1+a^{-1}+\dots+a^{-p+1}).$$

For the trefoil knot $T(2,3)$ the right hand  side of Equation (\ref{LuecksFormula}) yields:
\begin{eqnarray}
&&2\ln(k)-\sum_{n\geq 1} \frac 1 n \tr_{\C G} \left (1- \frac 1 {k^2} (1+a)(1+a^{-1})\right )^n
\label{Equation2} \\
&=& 2 \ln(k)-\sum_{n\geq 1} \frac 1 n \tr_{\C G} \left(\sum_{j=0}^n  {n \choose j} \left(- \frac 1 {k^2}\right )^j (2+ a+ a^{-1})^j  \right ) \nonumber
\end{eqnarray}

Hence, with a formal $\sqrt{a}$:

\begin{eqnarray*}
&=& 2 \ln(k)- \sum_{n\geq 1} \frac 1 n \tr_{\C G} \left(\sum_{j=0}^n {n \choose j} \left (- \frac 1 {k^2}\right ) ^j \left (\sqrt{a}+ \frac 1 {\sqrt{a}}\right )^{2 j}  \right )\\
&=& 2\ln(k) - \sum_{n\geq 1} \frac 1 n \sum_{j=0}^n  {n \choose j} \left (- \frac 1 {k^2}\right )^j {2 j \choose j} 
\end{eqnarray*}

On the other hand we know by Theorem \ref{MahlerMeasure} that Equation (\ref{LuecksFormula}) equals the
logarithm of the Mahler measure of $(1+a)$ which is $0$. This shows:

\begin{theorem} \label{LogarithmTheorem}
For $x  \geq 4$ a power series for the logarithm is given by:
$$\ln x = \sum_{n \geq 1} \sum_{j=0}^n \frac 1 n 
{n \choose j} {2 j \choose j}\left (-\frac 1 x \right )^j $$
\end{theorem}

\begin{remark}
Note, that Theorem \ref{LogarithmTheorem} does not require the powerful Theorem \ref{LuecksTheorem} but only
Theorem \ref{MahlerMeasure}. Thus, deeper topological arguments are not required.
\end{remark}

The case $x=4$ is somewhat special. Equation (\ref{Equation2}) at $x=k^2=4$ equals:

\begin{eqnarray*}
&&2\ln 2 -\sum_{n\geq 1} \frac 1 n \tr_{\C G} \left (\frac {a-2+a^{-1}}{-4}\right )^n\\
&=&2 \ln 2 - \sum_{n\geq 1} \frac 1 n \tr_{\C G} \left ( (-1)^n \frac {(\sqrt a - 1/\sqrt a)^{2 n}} {4^n}\right )\\
&=& 2 \ln 2 - \sum_{n \geq 1} \frac 1 n {2 n \choose n} \frac 1 {4^n}
\end{eqnarray*} 

Thus, we recover an identity which is well-known to Mathematica \cite{Mathematica} and was
derived by Lehmer in \cite{Lehmer:CentralBinomial}:

 
\begin{eqnarray}
\ln 4 &=& \sum_{n=1}^{\infty} \frac 1 n {2 n \choose n} \frac 1 {4^n}.
\end{eqnarray}
 
The function also converges at  
$x=k^2=2$ and with a similar computation we get:

$$
\ln 2 = \sum_{n=1}^{\infty} \frac 1 {2n} {2n \choose n} \frac 1 {4^n}.
$$

\subsection{A generating function}

It is interesting to look at the series for the natural logarithm in Theorem \ref{LogarithmTheorem} from a generating function point of view:

Let $$f_n(y)=\sum_{j=0}^n \frac 1 n {n \choose j} {2j \choose j} y^j$$
and $$F(x,y)=\sum_{n\geq 1} f_n(y) x^n$$ its generating function.

Similar to Example 5, Section 4.3 in \cite{Wilf:GeneratingFunctionology} we have:

\begin{eqnarray*}
F(x,y)&=&\sum_{n \geq 1} \sum_{j=1}^{n} \frac 1 n {n \choose j} {2 j \choose j} y^j x^n + \sum_{n \geq 1} \frac 1 n x^n\\
&=& \sum_{j=1}^{\infty} {2j \choose j} y^j \sum_{n=j}^{\infty} \frac 1 n {n \choose j} x^n - \ln (1-x)\\
&=& \sum_{j=1}^{\infty} {2j \choose j} y^j \frac 1 j \left ( \frac x {1-x} \right ) ^j - \ln (1-x)\\
&=& \sum_{j=1}^{\infty} {2j \choose j} \frac 1 j \left ( \frac {x y} {1-x} \right ) ^j - \ln (1-x)\\
&=& \ln 4 - 2 \ln \left (1 +\sqrt{1-4 \frac {x y} {1-x}}\right ) - \ln (1-x)\\
&=& \ln 4 - 2 \ln (\sqrt{1-x} + \sqrt{1-x- 4 xy})
\end{eqnarray*}

Here we make use of the identity \cite{Lehmer:CentralBinomial}
\begin{eqnarray*}
\sum_{j=1}^{\infty} \frac 1 j {2 j \choose j} z^j&=& 2 \log \left ( \frac {1-\sqrt{1-4 z}}{2 z} \right )
\end{eqnarray*}

\subsection{Triple Sum Series}

By working with a torus knot $T(p,q)$ and the matrix $A=(1+a+a^2+\dots+a^{p-1})$ in L\"uck's formula (\ref{LuecksFormula}) one readily obtains other infinite series for the logarithm.
For example if $p=3$ and thus $A=(1+a+a^2)$ we see for $k\geq \sqrt 3$:

\begin{eqnarray*}
2 \ln k &=& \sum_{n \geq 1} \frac 1 n \tr_{\C G} \left (1 - \frac 1 {k^2} (1+a+a^2) (1+a^{-1}+a^{-2})\right )^n\\
&=& \sum_{n \geq 1} \frac 1 n \tr_{\C G} \left ( \sum_{j=0}^n {n \choose j} \left ( - \frac 1 {k^2} \right )^j (a^{-1}+1+a)^{2 j}\right )\\
&=& \sum_{n \geq 1} \frac 1 n \sum_{j=0}^n {n \choose j} \left ( - \frac 1 {k^2} \right )^j \sum_{l=0}^j {2 j \choose {2 j - 2 l, l,l}},
\end{eqnarray*}
 
where the terms in the inner most sum are multinomial coefficients.

\section{Covering spaces}

For a given knot the terms in L\"ucks formula (\ref{LuecksFormula}) are by no means simple to compute. 
A single term depends on the chosen presentation. Furthermore, the convergence of the series is slow.              
An indirect approach via covering spaces of the knot complement gives more flexibility for computational simplifications. 
We illustrate it on the example $K$ being the trefoil knot.
It is well-known that the fundamental group of its knot complement
is isomorphic to the braid group $B_3$ on three strands:

$$B_3 = \langle \sigma_1, \sigma_2 | \sigma_1 \sigma_2 \sigma_1=\sigma_2 \sigma_1 \sigma_2 \rangle.$$

A direct application of Formula (\ref{LuecksFormula}) to this presentation would lead to a somewhat messy series.
However, the group $B_3$ has a natural homomorphic image in the symmetric group $S_3$ and its kernel is the pure braid group
$P_3$ of index $6$ in $B_3$.
By an application of Newworld's lemma \cite{DM:subgroup_separable} $P_3$ is isomorphic to the direct product of its center $C_3 \isom \Z$ and
a free group of rank $2$. More precisely:

With $t=(\sigma_1 \sigma_2)^3$ and $a=\sigma_1^2, b=\sigma_2^2$ a presentation for $P_3$ is given by

$$P_3 = \langle t,a,b | t a = a t, t b=b t \rangle.$$

Thus, by deleting the column corresponding to $t$ from the Fox Jacobian, the reduced matrix $A$ is:

$$A = \left ( \begin{array}{cc}
t-1 & 0\\
0 & t-1 
\end{array} \right ) .$$

By Theorem \ref{MahlerMeasure} the right hand-side of L\"ucks formula (\ref{LuecksFormula}) equals the logarithmic Mahler measure of
$\det(A)=(t-1)^2$, which is $0$. 

Now, since $P_3$ has index $6$ in $S^3-K$, for $K$ the trefoil, we know \cite{LueckBook:L2invariants} that 
Formula (\ref{LuecksFormula})  applied to the complement of the trefoil must be $\frac 1 6$ of this and thus also $0$.

\bigskip

\noindent
{\bf Acknowledgment: }  The author thanks Pat Gilmer, Walter Neumann and Neal Stoltzfus for valuable comments and suggestions.

\bibliography{../linklit}

\providecommand{\bysame}{\leavevmode\hbox to3em{\hrulefill}\thinspace}
\providecommand{\MR}{\relax\ifhmode\unskip\space\fi MR }
\providecommand{\MRhref}[2]{%
  \href{http://www.ams.org/mathscinet-getitem?mr=#1}{#2}
}
\providecommand{\href}[2]{#2}
\begin{thebibliography}{Wol99}

\bibitem[BZ85]{BZ}
G.~Burde and H.~Zieschang, \emph{Knots}, de Gruyter, Berlin, New York, 1985.

\bibitem[Den05]{Deninger:FugledeKadison}
Christopher Deninger, \emph{{Fuglede-Kadison determinants and entropy for
  actions of discrete amenable groups}}, 2005.

\bibitem[DM01]{DM:subgroup_separable}
Oliver~T. Dasbach and Brian~S. Mangum, \emph{The automorphism group of a free
  group is not subgroup separable}, Knots, braids, and mapping class
  groups---papers dedicated to Joan S. Birman (New York, 1998), AMS/IP Stud.
  Adv. Math., vol.~24, Amer. Math. Soc., Providence, RI, 2001, pp.~23--27.
  \MR{2002k:20063}

\bibitem[EW99]{EverestWard:MahlerMeasure}
Graham Everest and Thomas Ward, \emph{Heights of polynomials and entropy in
  algebraic dynamics}, Universitext, Springer-Verlag London Ltd., London, 1999.
  \MR{MR1700272 (2000e:11087)}

\bibitem[Leh85]{Lehmer:CentralBinomial}
D.~H. Lehmer, \emph{Interesting series involving the central binomial
  coefficient}, Amer. Math. Monthly \textbf{92} (1985), no.~7, 449--457.
  \MR{MR801217 (87c:40002)}

\bibitem[L{\"u}c94]{Lueck:Proceedings}
Wolfgang L{\"u}ck, \emph{{$L\sp 2$}-torsion and {$3$}-manifolds},
  Low-dimensional topology (Knoxville, TN, 1992), Conf. Proc. Lecture Notes
  Geom. Topology, III, Internat. Press, Cambridge, MA, 1994, pp.~75--107.
  \MR{MR1316175 (96g:57019)}

\bibitem[L{\"u}c02]{LueckBook:L2invariants}
\bysame, \emph{{$L\sp 2$}-invariants: theory and applications to geometry and
  {$K$}-theory}, Ergebnisse der Mathematik und ihrer Grenzgebiete. 3. Folge. A
  Series of Modern Surveys in Mathematics [Results in Mathematics and Related
  Areas. 3rd Series. A Series of Modern Surveys in Mathematics], vol.~44,
  Springer-Verlag, Berlin, 2002. \MR{1 926 649}

\bibitem[Wil94]{Wilf:GeneratingFunctionology}
Herbert~S. Wilf, \emph{generatingfunctionology}, second ed., Academic Press
  Inc., Boston, MA, 1994. \MR{MR1277813 (95a:05002)}

\bibitem[Wol99]{Mathematica}
Stephen Wolfram, \emph{The {M}athematica{$\sp \circledR$} book}, fourth ed.,
  Wolfram Media, Inc., Champaign, IL, 1999. \MR{MR1721106 (2000h:68001)}

\end{thebibliography}
\bibliographystyle {amsalpha}
\end{document}